## SPECIAL INVITED PAPER

## LARGE DEVIATIONS[1]

BY S. R. S. VARADHAN

*New York University*

This paper is based on Wald Lectures given at the annual meeting of the IMS in Minneapolis during August 2005. It is a survey of the theory of large deviations.

**1. Large deviations for sums.** The role of "large deviations" is best understood through an example. Suppose that $X_1, X_2, \ldots, X_n, \ldots$ is a sequence of i.i.d. random variables, for instance, normally distributed with mean zero and variance 1. Then,

$$E[e^{\theta(X_1+\cdots+X_n)}] = E[e^{\theta X_1}]^n = e^{n(\theta^2/2)}.$$

On the other hand,

$$E[e^{\theta(X_1+\cdots+X_n)}] = E[e^{n\theta(S_n/n)}].$$

Since, by the law of large numbers, $\frac{S_n}{n}$ is nearly zero, we have

$$E[e^{\theta(X_1+\cdots+X_n)}] = E[e^{o(n)}] \neq e^{o(n)}.$$

There is, of course, a very simple explanation for this. In computing expectations of random variables that can assume large values with small probabilities, contributions from such values cannot be ignored. After all, a product of something big and something small can still be big! In our case, assuming $\theta > 0$, for any $a > 0$,

$$E[e^{\theta S_n}] \geq e^{n\theta a} P\left[\frac{S_n}{n} \geq a\right] = e^{n\theta a} e^{-(na^2/2)+o(n)} = e^{n(a\theta-a^2/2)+o(n)}.$$

Received April 2007; revised April 2007.
[1]Supported in part by the NSF Grants DMS-01-04343 and DMS-06-04380.
*AMS 2000 subject classifications.* 60-02, 60F10.
*Key words and phrases.* Large deviations, rare events.

Since $a > 0$ is arbitrary,

$$E[e^{\theta S_n}] \geq e^{n \sup_{a>0}(\theta a - a^2/2) + o(n)} = e^{n\theta^2/2 + o(n)},$$

which is the correct answer.

The simplest example for which one can calculate probabilities of large deviations is coin tossing. The probability of $k$ heads in $n$ tosses of a fair coin is

$$P(n,k) = \binom{n}{k} 2^{-n} = \frac{n! 2^{-n}}{k!(n-k)!},$$

which using Stirling's approximation, is

$$\simeq \frac{\sqrt{2\pi} e^{-n} n^{n+1/2} 2^{-n}}{\sqrt{2\pi} e^{-(n-k)}(n-k)^{n-k+1/2} \sqrt{2\pi} e^{-k} k^{k+1/2}},$$

$$\log P(n,k) \simeq -\frac{1}{2}\log(2\pi)\left(n + \frac{1}{2}\right)\log n - \left(n - k + \frac{1}{2}\right)\log(n-k)$$

$$- \left(k + \frac{1}{2}\right)\log k - n\log 2$$

$$= -\frac{1}{2}\log(2\pi) - \frac{1}{2}\log n - \left(n - k + \frac{1}{2}\right)\log\left(1 - \frac{k}{n}\right)$$

$$- \left(k + \frac{1}{2}\right)\log\frac{k}{n} - n\log 2.$$

If $k \simeq nx$, then

$$\log P(n,k) \simeq -n[\log 2 + x\log x + (1-x)\log(1-x)] + o(n)$$
$$= -nH(x) + o(n),$$

where $H(x)$ is the Kullback–Leibler information or relative entropy of Binomial$(x, 1-x)$ with respect to Binomial$(\frac{1}{2}, \frac{1}{2})$.

This is not a coincidence. In fact, if $f_i$ are the observed frequencies in $n$ trials of a multinomial with probabilities $\{p_i\}$ for the individual cells, then

$$P(n, p_1, \ldots, p_k; f_1, \ldots, f_k) = \frac{n!}{f_1! \cdots f_k!} p_1^{f_1} \cdots p_k^{f_k}.$$

A similar calculation using Stirling's approximation yields, assuming $f_i \simeq nx_i$,

$$\log P(n, p_1, \ldots, p_k; f_1, \ldots, f_k) = -nH(x_1, \ldots, x_k; p_1, \ldots, p_k) + o(n),$$

where $H(x, p)$ is again the Kullback–Leibler information number

$$H(x) = \sum_{i=1}^{k} x_i \log \frac{x_i}{p_i}.$$

Any probability distribution can be approximated by one that is concentrated on a finite set and the empirical distribution from a sample of size $n$ will then have a multinomial distribution. One therefore expects that the probability $P(n, \alpha, \beta)$ that the empirical distribution

$$\frac{1}{n} \sum_{i=1}^{n} \delta_{X_i}$$

of $n$ independent observations from a distribution $\alpha$ is close to $\beta$ should satisfy

$$\log P(n, \alpha, \beta) = -nH(\beta, \alpha) + o(n),$$

where $H(\beta, \alpha)$ is again the Kullback–Leibler information number

$$\int \log \frac{d\beta}{d\alpha} \, d\beta = \int \frac{d\beta}{d\alpha} \log \frac{d\beta}{d\alpha} \, d\alpha.$$

This theorem, proven by Sanov, must be made precise. This requires a formal definition of what is meant by large deviations. We have a family $\{P_n\}$ of probability distributions on some space $X$ which we assume to be a complete separable metric space. There is a sequence of numbers $a_n \to \infty$ which we might as well assume to be $n$. Typically, $P_n$ concentrates around a point $x_0 \in X$ and, for sets $A$ away from $x_0$, $P_n(A)$ tends to zero exponentially rapidly in $n$, that is,

$$\log P_n(A) \simeq -nc(A),$$

where $c(A) > 0$ if $x_0 \notin A$.

We say that a *large deviation principle* holds for a sequence of probability measures $P_n$ defined on the Borel subsets of a Polish (complete separable metric) space $X$, with a *rate function* $H(x)$, if:

(1) $H(x) \geq 0$ is a lower semicontinuous function on $X$ with the property that $K_\ell = \{x : H(x) \leq \ell\}$ is a compact set for every $\ell < \infty$;

(2) for any closed set $C \subset X$,

$$\limsup_{n \to \infty} \frac{1}{n} \log P_n(C) \leq - \inf_{x \in C} H(x);$$

(3) for any open set $U \subset X$,

$$\liminf_{n \to \infty} \frac{1}{n} \log P_n(U) \geq - \inf_{x \in U} H(x).$$

While condition (1) is not really necessary for the validity of the large deviation principle, it is a useful condition on the rate function that will allow us to reduce the analysis to what happens on compact sets. Rate functions with this property are referred to as "good" rate functions.

For Sanov's theorem, the i.i.d. random variables will be taking values in a complete separable metric space $X$. Their common distribution will be a probability measure $\alpha \in \mathcal{M} = \mathcal{M}(X)$, the set of probability measures on $X$. The space $\mathcal{M}$ under weak convergence is a complete separable metric space. The empirical distributions will be random variables with values in $\mathcal{M}$, and $P_n$ will be their distribution, which will therefore be a probability measure on $\mathcal{M}$. The rate function

$$H(\beta, \alpha) = \int \frac{d\beta}{d\alpha} \log \frac{d\beta}{d\alpha} \, d\alpha$$

will be defined to be $+\infty$ unless $\beta \ll \alpha$ and the Radon–Nikodym derivative $f = \frac{d\beta}{d\alpha}$ is such that $f \log f$ is integrable with respect to $\alpha$.

**2. Rate functions, duality and generating functions.** We start with i.i.d. random variables and look at the sample mean

$$Y_n = \frac{S_n}{n} = \frac{X_1 + \cdots + X_n}{n}.$$

According to a theorem of Cramér, its distribution $P_n$ satisfies a large deviation principle with a rate function $h(a)$ given by

$$h(a) = \sup_\theta [\theta a - \log E[e^{\theta X}]].$$

In general, if a large deviation principle is valid for probability measures $P_n$ on a space $X$ with some rate function $H(x)$, it is not hard to see that, under suitable conditions on the function $F$ (boundedness and continuity will suffice),

$$\frac{1}{n} \log \int e^{nF(x)} \, dP_n \to \sup_x [F(x) - H(x)].$$

The basic idea is just

$$\frac{1}{n} \log \sum e^{na_i} \to \sup_i a_i.$$

In other words, the logarithms of generating functions are dual to, or Legendre transforms of, large deviation rate functions. For instance,

$$\log \left[ \int e^{V(x)} \, d\alpha \right] = \sup_\beta \left[ \int V(x) \, d\beta - H(\beta, \alpha) \right].$$

If the rate function is convex, as it is in the case of sums of independent identically distributed random variables, the duality relationship is invertible and

$$H(\beta, \alpha) = \sup_{V(\cdot)} \left[ \int V(x) \, d\beta - \log E[e^{V(x)} \, d\alpha] \right],$$

where the supremum on $V$ in either case is taken over all bounded measurable functions or bounded continuous functions.

**3. Markov processes.** This relationship for i.i.d. sequences can be extended to the Markovian context. For simplicity, let us assume that we have a finite state space $X$ and transition probabilities $\pi(x,y)$ of a Markov chain on $X$. Let us suppose that $\pi(x,y) > 0$ for all $x, y$. If $V(\cdot): X \to R$ is a function on $X$, then

$$E_x[\exp[V(X_1) + V(X_2) + \cdots + V(X_n)]]$$

can be explicitly evaluated as

$$\sum_y \pi_V^n(x,y),$$

where $\pi_V(x,y) = \pi(x,y)e^{V(y)}$ and $\pi_V^n$ is the $n$th power of $\pi_V$. Since $\pi_V$ is a matrix with positive entries,

$$\frac{1}{n} \log \sum_y \pi_V^n(x,y) \to \log \lambda_\pi(V),$$

where $\lambda_\pi(V)$ is the principal eigenvalue of $\pi_V$. The analog of Sanov's theorem in this context establishes a large deviation result for the empirical distributions of a long sequence $\{x_1, x_2, \ldots, x_n\}$. They belong to the space $\mathcal{M}$ of probability measures on $X$. Let us denote by $P$ the Markov process with transition probability $\pi(x,y)$ and by $Q_n$ the measure on $\mathcal{M}$ which is the distribution of the empirical distribution $\{x_1, x_2, \ldots, x_n\}$. The large deviation upper bound for $Q_n$, with a rate function given by

$$H_\pi(q) = \sup_V \left[ \sum_x V(x)q(x) - \log \lambda_\pi(V) \right],$$

is an easy consequence of estimates on the generating function.

There is a more direct way of approaching $H_\pi$ through lower bounds. Let us pretend that our Markov chain exhibits "atypical" behavior and behaves like a different Markov chain, one with transition probability $\tilde{\pi}(x,y)$. In other words, the empirical distribution of visits to different sites of $X$ is close to $q$, which turns out to be the invariant distribution for a different chain, one with transition probabilities $\tilde{\pi}(x,y)$. We can estimate the probability of this event. Let $U$ be an open set around $q$ in the space $\mathcal{M}$ of probability measures on $X$. Let $A_n$ be the set of all realizations of $\{x_1, \ldots, x_n\}$ with empirical distributions belonging to $U$. We can estimate $P(A_n)$, the probability of $A_n$ under the original $\pi(x,y)$ chain, by

$$P(A_n) = \sum_{x_1, x_2, \ldots, x_n \in A_n} \pi(x, x_1) \cdots \pi(x_{n-1}, x_n)$$

$$= \sum_{x_1, x_2, \ldots, x_n \in A_n} \tilde{\pi}(x, x_1) \cdots \tilde{\pi}(x_{n-1}, x_n) \exp\left[ -\sum_{i=1}^n \log \frac{\tilde{\pi}(x_{i-1}, x_i)}{\pi(x_{i-1}, x_i)} \right]$$

$$= \int_{A_n} \exp\left[-\sum_{i=1}^n \log \frac{\tilde\pi(x_{i-1},x_i)}{\pi(x_{i-1},x_i)}\right] d\tilde P,$$

where $\tilde P$ is the Markov chain with transition probability $\tilde\pi(x,y)$. By the ergodic theorem,

$$\lim_{n\to\infty} \tilde P(A_n) = 1$$

as $n \to \infty$. An application of Jensen's inequality yields

$$\liminf_{n\to\infty} \frac{1}{n} \log Q_n(A) \geq -\sum_{x,y} \tilde\pi(x,y) \log \frac{\tilde\pi(x,y)}{\pi(x,y)} q(x).$$

We can pick any $\tilde\pi$, provided $q$ is invariant for $\tilde\pi$, that is, $q\tilde\pi = q$ and we will get an upper bound for $H_\pi(q)$. Therefore,

$$H_\pi(q) \leq \inf_{\tilde\pi:\tilde q=q} \sum_{x,y} \tilde\pi(x,y) \log \frac{\tilde\pi(x,y)}{\pi(x,y)} q(x).$$

In fact, there is equality here and $H_\pi(q)$ is dual to $\log \lambda_\pi(V)$.

$$H_\pi(q) = \sup_{V(\cdot)} \left[\sum_x V(x)q(x) - \log \lambda_\pi(V)\right],$$

$$\log \lambda_\pi(V) = \sup_{q(\cdot)} \left[\sum_x V(x)q(x) - H_\pi(q)\right].$$

If we are dealing with a process in continuous time, we will have a matrix $A$ of transition rates $\{a(x,y)\}$ with $a(x,y) \geq 0$ for $x \neq y$ and $\sum_y a(x,y) = 0$. The transition probabilities $\{\pi(t,x,y)\}$ will be given by $\pi(t,x,y) = (\exp tA)(x,y)$. If $a$ is symmetric, that is, $a(x,y) = a(y,x)$, then so is $\pi(t,\cdot,\cdot)$. The uniform distribution on $X$ will, in this case, be the invariant measure. The role of $\log \lambda(V)$ in the previous discrete situation will now be played by the principal eigenvalue $\lambda_a(V)$ of $A+V$, where $(A+V)(x,y) = a(x,y) + \delta(x,y)V(y)$. Here, $I = \{\delta(x,y)\}$ is the identity matrix. This is the conjugate of the rate function $H(q)$.

$$H_a(q) = \sup_{V(\cdot)} \left[\sum_x V(x)q(x) - \lambda_a(V)\right],$$

$$\lambda_a(V) = \sup_{q(\cdot)} \left[\sum_x V(x)q(x) - H_a(q)\right].$$

Note that, in the symmetric case, we have the usual variational formula

$$\lambda_a(V) = \sup_{\substack{u(\cdot) \\ \sum_x [u(x)]^2=1}} \left[\sum_x V(x)u^2(x) + \sum_{x,y} a(x,y)u(x)u(y)\right]$$

$$= \sup_{\substack{u(\cdot) \\ \sum_x [u(x)]^2 = 1}} \left[ \sum_x V(x) u^2(x) - \tfrac{1}{2} \sum_{x,y} a(x,y)(u(x) - u(y))^2 \right]$$

$$= \sup_{\substack{q(\cdot) \geq 0 \\ \sum_x q(x) = 1}} \left[ \sum_x V(x) u^2(x) - \tfrac{1}{2} \sum_{x,y} a(x,y)(\sqrt{q(x)} - \sqrt{q(y)})^2 \right].$$

It is therefore not very surprising that

$$H_a(q) = \tfrac{1}{2} \sum_{x,y} a(x,y)(\sqrt{q(x)} - \sqrt{q(y)})^2.$$

**4. Small random perturbations. The exit problem.** The context for large deviations is a probability distribution that is nearly degenerate. So far, this has come from the law of large numbers or the ergodic theorem. But it can also come from small random perturbations of a deterministic system.

Consider, for instance, the Brownian motion or Wiener measure on the space of continuous functions $C_0[[0,T];R]$ that are zero at 0, that is, $x(0) = 0$. We can make the variance of the Brownian motion at time $t$ equal to $\epsilon t$ instead of $t$. We can consider $x_\epsilon(t) = \sqrt{\epsilon} x(t)$. Perhaps consider $\epsilon = \frac{1}{n}$ and $x_n(t) = \frac{1}{n} \sum_i y_i(t)$, the average of $n$ independent Brownian motions. In any case, the measure $P_\epsilon$ of $x_\epsilon$ is nearly degenerate at the path $f(\cdot) \equiv 0$.

We have a large deviation principle for $P_\epsilon$ with a rate function

$$H(f) = \tfrac{1}{2} \int_0^T [f'(t)]^2 \, dt.$$

There are various ways of seeing this. We will exhibit two. As a Gaussian process, the generating function for Brownian motion is

$$\log E \left[ \exp \left[ \int_0^T x(t) g(t) \, dt \right] \right] = \tfrac{1}{2} \int_0^T \int_0^T \min(s,t) g(s) g(t) \, ds \, dt.$$

Its dual is given by

$$\sup_{g(\cdot)} \left[ \int_0^T f(t) g(t) \, dt - \tfrac{1}{2} \int_0^T \int_0^T \min(s,t) g(s) g(t) \, ds \, dt \right]$$

$$= \tfrac{1}{2} \int_0^T [f'(t)]^2 \, dt.$$

Or, we can perturb $\sqrt{\epsilon} x(t)$ by $f(t)$ and change the measure from $P_\epsilon$ to $Q_\epsilon$ that concentrates near $f$ rather that at 0. Of course, there are many measures that do this. $Q_\epsilon$ is just one such. The relative entropy is easily calculated for our choice. It is

$$H(Q_\epsilon, P_\epsilon) = \frac{1}{2\epsilon} \int_0^T [f'(t)]^2 \, dt.$$

This is also a lower bound for the rate function $H(f)$. In general, the rate function has a lower bound which is the "entropy cost" of changing the measure to do what we want it to do, namely, to concentrate on the "wrong spot." The cheapest way of achieving this is invariably the actual rate function. Let us now start with an ODE in $R^n$,

$$dx(t) = b(x(t))\,dt, \qquad x(0) = x_0,$$

and perturb it with a small noise,

$$dx_\epsilon(t) = b(x_\epsilon(t))\,dt + \sqrt{\epsilon}\beta(t), \qquad x(0) = x_0,$$

where $\beta(t)$ is the Brownian motion in $R^n$. As $\epsilon \to 0$, the distribution $P_\epsilon$ of $x_\epsilon(\cdot)$ will concentrate on the unique solution $x_0(\cdot)$ of the ODE. But there will be a large deviation principle with a rate function

$$H(f) = \tfrac{1}{2} \int_0^T \|f'(t) - b(f(t))\|^2\,dt$$

which will happen if $\sqrt{\epsilon}\beta(\cdot)$ concentrates around $g(\cdot)$ with $g'(\cdot) = f'(\cdot) - b(f(\cdot))$.

One application of this is to the "exit problem" and the resulting interpretation of "punctuated equilibria." Let us take the ODE to be a gradient flow:

$$dx(t) = -(\nabla V)(x(t))\,dt.$$

The system will move toward a minimum of $V$. If $V$ has multiple local minimum or valleys surrounded by mountains, the solutions of the ODE could be trapped near a local minimum, depending on the starting point. They will not move from one local minimum to a deeper minimum. On the other hand, with even a small noise, they develop "wanderlust." While, most of the time, they follow the path dictated by the ODE, they will, from time to time, deviate sufficiently to be able to find and reach a lower minimum. If one were to sit at a deeper minimum and wait for the path to arrive, then trace its history back, one would be apt to find that it did not wander at all, but made the most efficient beeline from the previous minimum to this one, as if it were guided by a higher power. The large deviation interpretation of this curious phenomenon is simply that the system will experiment with every conceivable path, with probabilities that are extremely small. The higher the rate function, that is, the less efficient the path, the smaller will be the probability, hence, less frequent the attempts involving that path. Therefore, the first attempt to take the path to its new location will be the most efficient path, one that climbs the lowest mountain pass to get to the new fertile land!

**5. Gibbs measures and statistical mechanics.** Let us look at point processes on $R$. The simplest is the Poisson point process. One can imagine it as starting from a finite interval $[-\frac{n}{2}, \frac{n}{2}]$ and placing $k_n = \rho n$ particles there randomly, independently and uniformly. Their joint density is

$$\frac{1}{n^{k_n}} dx_1 \cdots dx_{k_n}.$$

We can ignore the labeling and think of it as a point process $P_n$. The density is then

$$\frac{k_n!}{n^{k_n}} dx_1 \cdots dx_{k_n}.$$

As $n \to \infty$, we obtain a point process $P$ which is the Poisson point process with intensity $\rho$. We can try to modify $P_n$ into $Q_n$ given by

$$dQ_n = \frac{1}{Z_n} \exp\left[-\sum_{i,j} V(x_i - x_j)\right] dP_n,$$

where $V \geq 0$ has compact support. Here, $Z_n$ is the normalizing constant

$$Z_n = \int \exp\left[-\sum_{i,j} V(x_i - x_j)\right] dP_n.$$

The limit

$$\lim_{n \to \infty} \frac{1}{n} \log Z_n = \Psi(V)$$

is called the *free energy*. If $Q$ is any stationary, that is, translation-invariant point process on $R$, its restriction to $[a, b]$ is a symmetric (i.e., permutation-invariant) probability distribution $q_{[a,b]}$ on

$$\bigcup_{k \geq 0} [a, b]^k.$$

Similarly, the Poisson point process with intensity $\rho$ generates $p_{[a,b]}$, a convex combination of uniform distributions on $[a, b]^k$ with Poisson weights

$$e^{-\rho(b-a)} \frac{\rho^k (b-a)^k}{k!}.$$

The relative entropy $H(q_{[a,b]}, p_{[a,b]}) = H_\ell(Q, P)$ depends only on the length $\ell = (b - a)$ and is superadditive, that is,

$$H_{\ell_1 + \ell_2}(Q, P) \geq H_{\ell_1}(Q, P) + H_{\ell_2}(Q, P).$$

Therefore, the limit

$$\lim_{\ell \to \infty} \frac{1}{\ell} H_\ell(Q, P) = \sup_\ell H_\ell(Q, P) = \widehat{H}(Q, P)$$

exists and is called the *specific relative entropy*. There is the expectation

$$E^Q[V] = \lim_{n\to\infty} \frac{1}{n} E^Q\left[\sum_{x_\alpha, x_\beta \in [-n/2, n/2]} V(x_\alpha - x_\beta)\right].$$

Then,

$$\Psi(V) = -\inf_Q [E^Q[V] + \widehat{H}(Q,P)]$$

over all stationary point processes and the minimizer (unique) is the limit of $Q_n$. Such considerations play a crucial role in the theory of equilibrium statistical mechanics and thermodynamics [16].

**6. Interacting particle systems.** Interacting particle systems offer an interesting area where methods of large deviation can be applied. Let us look at some examples. Consider the lattice $\mathbf{Z}_N$ of integers modulo $N$. We think of them as points arranged uniformly on the circle of unit circumference or the unit interval with endpoints identified. For each $x \in \mathbf{Z}_N$, there corresponds a point $\xi = \frac{x}{N}$ on the interval (really, the circle). There are a certain number of particles in the system that occupy some of the sites. The particles wait for an exponential time and then decide to jump. They pick the adjacent site either to its left or right with equal probability and then jump. In this model, the particles do not interact and, after a diffusive rescaling of space and time, the particles behave like independent Brownian motions. Since there is a large number of particles, by an application of the law of large numbers, the empirical density will be close to a solution $\rho(t, \xi)$ of the heat equation

$$\rho_t = \tfrac{1}{2}\rho_{\xi\xi},$$

with a suitable initial condition, depending on how the particles were initially placed.

We can introduce an interaction in this model, by imposing a limit of at most one particle per site. When the particle decides to jump to a randomly chosen site after a random exponential waiting time, it can jump only when the site is unoccupied. Otherwise, it must wait for a new exponential time before trying to jump again. All of the particles wait and attempt to jump independently of each other. Since we are dealing with continuous-time and exponential distributions, two particle will not try to jump at the same time and we do not have to resolve ties.

One can express all of this by simply writing down the generator. $\eta$ is a configuration and $X_N = \{0,1\}^N$ is the set of all configurations. $\eta(x) = 1$ if the site $x$ is occupied and $\eta(x) = 0$ if the site $x$ is free. A point in $X_N$ is just a map $\eta : \mathbf{Z}_N \to \{0,1\}$. The generator $A_N$ acting on functions $f : X_N \to R$ is given by

$$(A_N f)(\eta) = \tfrac{1}{2} \sum_x \eta(x)[(1 - \eta(x+1))(f(\eta^{x,x+1}) - f(\eta))$$
$$+ (1 - \eta(x-1))(f(\eta^{x,x-1}) - f(\eta))],$$

where $\eta^{x,y}$ is defined by

$$\eta^{x,y}(z) = \begin{cases} \eta(x), & \text{if } z = y, \\ \eta(y), & \text{if } z = x, \\ \eta(z), & \text{otherwise.} \end{cases}$$

This system has some interesting features. The total number $\sum_x \eta(x)$ of particles does not change over time. Particles just diffuse over time. It will take time of order $N^2$ for the effect to be felt at a distance of order $N$. The equilibrium distributions $\mu_{N,k}$ are uniform over all possible $\binom{N}{k}$ configurations, where $k$ is the number of particles. In particular, there are multiple equilibria. Such systems can be locally in equilibrium while approaching the global equilibrium rather slowly. For example, we initially place our particles in such a way that one half of the circle $C_1$ has particle density $\tfrac{1}{4}$, while the other half $C_2$ has density $\tfrac{3}{4}$. The system will locally stay near the two different equilibria at the two intervals for quite some time, with just some smoothing near the edges. If we wait for times of order $N^2$, that is, time $N^2 t$ with $t > 0$, the density will be close to $\rho(t, \tfrac{x}{N})$. A calculation that involves applying the speeded up generator $N^2 A_N$ to expressions of the form

$$\frac{1}{N} \sum J\left(\frac{x}{N}\right) \eta(x)$$

shows that $\rho(t, \xi)$ will be the solution of the heat equation

$$\rho_t(t, \xi) = \tfrac{1}{2} \rho_{\xi\xi}(t, \xi); \qquad \rho(0, \xi) = \tfrac{1}{4} \mathbf{1}_{C_1}(\xi) + \tfrac{3}{4} \mathbf{1}_{C_2}(\xi).$$

In this case, this behavior is the same as if the particles moved independently without any interaction. While this is true at the level of the law of large numbers, we shall see that for large deviation rates, the interaction does matter.

**7. Superexponential estimates.** Let us consider the process in equilibrium, with $k = \rho N$ particles. We let $N, k$ go to $\infty$ while $\rho$ is kept fixed. It is easy to see, for any fixed site $x$, that

$$\mu_{N,k}[\eta(x) = 1] \to \rho$$

and $\{\eta(x_j)\}$ become independent for any finite set of distinct sites. If $f(\eta)$ is a local function and $f_x(\eta) = f(\tau_x \eta)$ its translation, then

$$\frac{1}{N} \sum_x f_x(\eta) \to E^{P_\rho}[f(\eta)]$$

as $N \to \infty$. The probabilities of deviations decay exponentially fast. The exponential rate is given by

$$H(a) = \inf_{Q \in \mathcal{M}(a,\rho)} H(Q, P_\rho),$$

where $\mathcal{M}(a, \rho)$ consists of stationary measures $Q$ with density $\rho$ and $E^Q[f(\eta)] = a$, and $P_\rho$ is Bernoulli with density $\rho$. $H(Q, P)$ is the specific entropy, calculated as the limit

$$\lim_{n \to \infty} \frac{1}{n} H(Q_n, P_n),$$

where $Q_n$ and $P_n$ are restrictions of $Q$ and $P$ to a block of $n$ sites. This is, of course, an equilibrium calculation done at a fixed time. What about space–time averages

$$\frac{1}{NT} \int_0^T \sum_x f_x(\eta(t)) \, dt?$$

The rate function is now some $H_T(a)$ that depends on $T$. While this is hard to calculate, one can show that

$$\lim_{T \to \infty} H_T(a) = \inf_{\rho(\cdot) \in \mathcal{P}(a,\rho)} H(P_{\rho(\xi)}, P_\rho) \, d\xi,$$

where $\mathcal{P}(a, \rho)$ consists of density profiles $\rho(\xi)$ that satisfy

$$\int \rho(\xi) \, d\xi = \rho \quad \text{and} \quad \int E^{P_{\rho(\xi)}}[f(\eta)] \, d\xi = a.$$

Note that the rate increases to a finite limit as $T \to \infty$. The reason is that large-scale fluctuations in density can occur with exponentially small probability and these fluctuations do not necessarily decay when $T$ is large. Particles diffuse slowly and it takes times of order $N^2$ to diffuse through order $N$ sites.

To resolve this difficulty, we define the approximate large-scale empirical density by

$$r_{N,\epsilon}(\xi, \eta) = \frac{1}{N\epsilon} \sum_{x: |x/N - \xi| \leq \epsilon} \eta(x)$$

and compare

$$\frac{1}{N^3} \int_0^{N^2 T} \sum_x f_x(\eta(t)) \, dt = Y_{N,f}(T)$$

with

$$\frac{1}{N^3} \int_0^{N^2 T} \sum_x \widehat{f}\left(r_{N,\epsilon}\left(\frac{x}{N}, \eta(t)\right)\right) dt = \widehat{Y}_{N,f,\epsilon}(T),$$

where
$$\widehat{f}(\rho) = E^{P_\rho}[f(\eta)].$$

One can show that
$$\limsup_{\epsilon \to 0} \limsup_{N \to \infty} \frac{1}{N} \log P[|Y_{N,f}(T) - \widehat{Y}_{N,f,\epsilon}(T)| \geq \delta] = -\infty$$

for any $\delta > 0$. In other words, in time scale $N^2$, the large scale density fluctuations are responsible for all large deviations with a logarithmic decay rate of $N$. Such estimates, called *superexponential estimates*, play a crucial role in the study of interacting particle systems where approach to equilibrium is slow in large scales.

**8. Hydrodynamical limits.** If we consider a large system with multiple equilibria, the system can remain in a local equilibrium for a long time. Let the different equilibria be labeled by parameter $\rho$. If the system lives on a large spatial domain and if $\xi$ denotes a point in the macroscopic scale of space, a function $\rho(\xi)$ could describe a system that is locally in equilibrium, albeit at different ones at macroscopically different points. If $t$ is time measured in a suitably chosen faster scale, then $\rho(\xi) = \rho(t, \xi)$ may evolve gently in time and converge as $t \to \infty$ to a constant $\rho$ that identifies the global equilibrium.

For instance, in our example, let us initially start with $k = N\rho$ particles and distribute them in such a way that we achieve a density profile $\rho_0(\xi)$ as $N \to \infty$. Technically, this means
$$\lim_{N \to \infty} \frac{1}{N} \sum_x J\left(\frac{x}{N}\right) \eta(x) = \int J(\xi) \rho_0(\xi) \, d\xi$$

for bounded continuous test functions $J$. What will we see after time $tN^2$, especially when $N$ is large?

The answer is very easy. Consider the sum
$$f(\eta) = \frac{1}{N} \sum_x J\left(\frac{x}{N}\right) \eta(x)$$

and apply the speeded up operator $N^2 A_N$ to it.
$$f(\eta^{x,y}) - f(\eta) = \frac{1}{N}[\eta(x) - \eta(y)]\left[J\left(\frac{y}{N}\right) - J\left(\frac{x}{N}\right)\right],$$
$$(N^2 A_N f)(\eta) = \frac{N}{2} \sum_x \eta(x)(1 - \eta(x+1))\left[J\left(\frac{x+1}{N}\right) - J\left(\frac{x}{N}\right)\right]$$
$$+ \eta(x)(1 - \eta(x-1))\left[J\left(\frac{x-1}{N}\right) - J\left(\frac{x}{N}\right)\right]$$

$$= \frac{N}{2} \sum_x [\eta(x)(1-\eta(x+1)) - \eta(x+1)(1-\eta(x))]$$
$$\times \left(J\left(\frac{x+1}{N}\right) - J\left(\frac{x}{N}\right)\right)$$
$$= \frac{N}{2} \sum_x [\eta(x) - \eta(x+1)]\left[J\left(\frac{x+1}{N}\right) - J\left(\frac{x}{N}\right)\right]$$
$$= \frac{N}{2} \sum_x \left[J\left(\frac{x+1}{N}\right) + J\left(\frac{x-1}{N}\right) - 2J\left(\frac{x}{N}\right)\right]\eta(x)$$
$$\simeq \frac{1}{2N} \sum_x J''\left(\frac{x}{N}\right)\eta(x),$$

leading to the limiting heat equation

$$\rho_t = \tfrac{1}{2}\rho_{\xi\xi}, \qquad \rho(0,\xi) = \rho_0(\xi).$$

Let us change the problem slightly. Introduce a slight bias. The probabilities to right and left are $\frac{1}{2} \pm \frac{1}{2N}b$, where $b > 0$ is the bias to the right. This introduces an extra term, so that now

$$N^2 A_N \frac{1}{N} \sum_x J\left(\frac{x}{N}\right)\eta(x) = \frac{1}{2N} \sum_x J''\left(\frac{x}{N}\right)\eta(x) + F_N(\eta),$$

where

$$F_N(\eta) = \frac{b}{2N} \sum J'\left(\frac{x}{N}\right)[\eta(x)(1-\eta(x+1)) + \eta(x)(1-\eta(x-1))].$$

This is a problem because $F_N$ is nonlinear in $\eta$ and cannot be replaced by a simple expression involving $\rho$. If we were locally in equilibrium, $F_N$ would be replaced by

$$G_{N,\epsilon}(\eta) = \frac{b}{N} \sum J'\left(\frac{x}{N}\right) r_{N,\epsilon}\left(\frac{x}{N},\eta\right)\left(1 - r_{N,\epsilon}\left(\frac{x}{N},\eta\right)\right),$$

where $r_{N,\epsilon}(\frac{x}{N},\eta)$ is the empirical density. If we were in a global equilibrium, in the faster time scale,

$$\left|\int_0^T F_N(\eta(t)) - G_{N,\epsilon}(\eta(t))\,dt\right| \geq \delta,$$

with superexponentially small probability. The importance of superexponential estimates lies in the following elementary, but universal, inequality:

$$Q(A) \leq \frac{2 + H(Q,P)}{\log(1/P(A))}.$$

If $P(A)$ is superexponentially small and $H(Q,P)$ is linear in $N$, then $Q(A)$ is small. Perturbation of the process by a slight bias of order $\frac{1}{N}$ produces a relative entropy of order $\frac{1}{N^2}$ per site per unit time. With $N$ sites and time of order $N^2$, this is still only linear in $N$. Changes in initial conditions can also only contribute linear relative entropy. So, we can carry out the approximation and end up with

$$\rho_t = \tfrac{1}{2}\rho_{\xi\xi} - [b\rho(1-\rho)]_\xi.$$

Note that without exclusion, the limit would have been

$$\rho_t = \tfrac{1}{2}\rho_{\xi\xi} - [b\rho]_\xi,$$

which is just the CLT when there is a small mean.

**9. Large deviations in hydrodynamical limits.** We saw earlier that in our model, if we started with an initial profile with density $\rho_0$, in a speeded up time scale, the system will evolve with a time-dependent profile which is the solution of the heat equation

$$\rho_t = \tfrac{1}{2}\rho_{\xi\xi}.$$

This is, of course, true with probability nearly 1. We can make the initial condition deterministic in order that no initial deviation is possible. Still, our random evolution could produce, with small probability, strange behavior. If we limit ourselves to deviations with only exponentially small probabilities, what deviations are possible and what are their exponential rates?

We can doctor the system with a bias $b$ which is not constant, but is a function $b(t, \frac{x}{N})$ of $x$ and fast time $t$ with $b(t,\xi)$ being a nice function of $t$ and $\xi$. This will produce a solution of

$$\rho_t(t,\xi) = \tfrac{1}{2}D_\xi^2 \rho(t,\xi) - D_\xi[b(t,\xi)\rho(t,\xi)(1-\rho(t,\xi))], \qquad \rho(0,\xi) = \rho_0(\xi).$$

This can be done with an entropy cost that can be calculated as

$$\Psi(b) = \tfrac{1}{2}\int_0^T \int |b(t,\xi)|^2 \rho(t,\xi)(1-\rho(t,\xi))\, dt\, d\xi,$$

for the duration $[0,T]$. If $\rho(t,\xi)$ is given, then we can minimize $\Psi(b)$ over all compatible $b$. This is a lower bound for the rate. One can match this bound in the other direction.

**10. Large deviations for random walks in random environments.** We will start with a probability space $(\Omega, \Sigma, P)$ on which $\mathbf{Z}^d$ acts ergodically as a family $\tau_z$ of measure-preserving transformations. We are given $\pi(\omega, z)$, which is a probability distribution on $\mathbf{Z}^d$ for each $\omega$, and is a measurable functions of $\omega$ for each $z$. One can then generate random transition probabilities $\pi(\omega, z', z)$ by defining

$$\pi(\omega, z', z' + z) = \pi(\tau_{z'}\omega, z).$$

For each $\omega$, $\pi(\omega, z', z)$ can serve as the transition probability of a Markov process on $\mathbf{Z}^d$, and the measure corresponding to this process, starting from 0, is denoted by $Q^\omega$. This is, of course, random (it depends on $\omega$ and is called the *random walk in the random environment* $\omega$). One can ask the usual questions about this random walk and, in some form, they may be true for almost all $\omega$ with respect to $P$. The law of large numbers, if valid, will take the form

$$P\left[\omega : \lim_{n\to\infty} \frac{S_n}{n} = m(P) \text{ a.e. } Q^\omega\right] = 1.$$

Such statement concerning the almost sure behavior under $Q^\omega$ for almost all $\omega$ with respect to $P$ are said to deal with the "quenched" version. Sometimes one wishes to study the behavior of the "averaged" or "annealed" measure,

$$\overline{Q} = \int Q^\omega P(d\omega).$$

The law of large numbers is the same because it is equivalent to

$$\overline{Q}\left[\omega : \lim_{n\to\infty} \frac{S_n}{n} = m(P)\right] = 1.$$

On the other hand, questions on the asymptotic behavior of probabilities, like the central limit theorem or large deviations, could be different for the quenched and the averaged cases.

A special environment, which is called the *product environment*, is one in which $\pi(\omega, z', z' + z)$ are independent for different $z'$ and have a common distribution $\beta$ which is a probability measure on the space $\mathcal{M}$ of all probability measures on $\mathbf{Z}^d$. In this case, the canonical choice for $(\Omega, \Sigma, P)$ is the countable product of $\mathcal{M}$ and the product measure $P$ with marginals $\beta$.

There are large deviation results regarding the limits

$$\lim_{n\to\infty} \frac{1}{n} \log Q^\omega\left[\frac{S_n}{n} \simeq a\right] = I(a)$$

and

$$\lim_{n\to\infty} \frac{1}{n} \log \overline{Q}\left[\frac{S_n}{n} \simeq a\right] = \overline{I}(a).$$

The difference between $I$ and $\overline{I}$ has a natural explanation. In the special case of $Z$, in terms of the large deviation behavior in the one-dimensional random environment, it is related to the following question. If, at time $n$, we see a particle reach an improbable value $na$, what does it say about the environment in $[0, na]$? Did the particle behave strangely in the environment or did it encounter a strange environment? It is probably a combination of both.

Large deviation results exists, however, in much wider generality, both in the quenched and the averaged cases. The large deviation principle is essentially the existence of the limits

$$\lim_{n\to\infty} \frac{1}{n} \log E[\exp[\langle \theta, S_n \rangle]] = \Psi(\theta).$$

The expectation is with respect to $Q^\omega$ or $\overline{Q}$, which could produce different limits for $\Psi$. The law of large numbers and central limit theorem involve the differentiability of $\Psi$ at $\theta = 0$, which is harder. In fact, large deviation results have been proven by general subadditivity arguments for the quenched case. Roughly speaking, fixing $\omega$, we have

$$Q^\omega[S_{k+\ell} \simeq (k+\ell)a] \geq Q^\omega[S_k \simeq ka] \times Q^{\tau_{ka}\omega}[S_\ell \simeq \ell a].$$

It is harder for the averaged case. We want to prove that the limits

$$\lim_{n\to\infty} \frac{1}{n} \log \overline{Q}\left[\frac{S_n}{n} \simeq a\right] = -\overline{I}(a)$$

or, equivalently,

$$\lim_{n\to\infty} \frac{1}{n} \log E^{\overline{Q}}[\exp[\langle \theta, S_n \rangle]] = \overline{\Psi}(\theta),$$

exist. The problem is that the measure $\overline{Q}$ is not very nice. As the random walk explores $\mathbf{Z}^d$, it learns about the environment and in the case of the product environment, when it returns to a site that it has visited before, the experience has not been forgotten and leads to long term correlations. However, if we are interested in the behavior $S_n \simeq na$ with $a \neq 0$, the same site is not visited too often and the correlations should rapidly decay.

One can use Bayes' rule to calculate the conditional distribution

$$\overline{Q}[S_{n+1} = S_n + z | S_1, S_2, \ldots, S_n] = q(z|w),$$

where $w$ is the past history of the walk. Before we do that, it is more convenient to shift the origin as we go along so that the current position of the random walk is always the origin and the current time is always 0. An $n$ step walk then looks like $w = \{S_0 = 0, S_{-1}, \ldots, S_{-n}\}$. We pick a $z$ with probability $q(z|w)$. We obtain a new walk of $n+1$ steps $w' = \{S'_0 = 0, S'_{-1}, \ldots, S'_{-(n+1)}\}$ given by $S'_{-(k+1)} = S_{-k} - z$ for $k \geq 0$. We can now calculate $q(z|w)$. We need to know all the quantities $\{k(w,x,z)\}$, the numbers of times the walk has visited $x$ in the past and jumped from $x$ to $x+z$. It is not hard to see that the a posteriori probability can be calculated as

$$q(z|w) = \frac{\int \pi(z) \Pi_{z'} \pi(z')^{k(w,0,z')} \beta(d\pi)}{\int \Pi_{z'} \pi(z')^{k(w,0,z')} \beta(d\pi)}.$$

While this initially makes sense only for walks of finite length, it can clearly be extended to all transient paths. Note that although we only use $k(w, 0, z')$, in order to obtain the new $k(w', 0, z')$, we would need to know the collection $\{k(w, z, z')\}$.

Now, suppose that $R$ is a process with stationary increments $\{z_j\}$. Then, we can again make the current position the origin and if the process is transient, as it would be if the increments were ergodic and had a nonzero mean $a$, the conditional probabilities $q(z|w)$ would exist a.e. $R$ and can be compared to the corresponding conditional probabilities $r(z|w)$ under $R$. The relative entropy

$$H(R) = E^R\left[\sum_z r(z|w) \log \frac{r(z|w)}{q(z|w)}\right]$$

is then well defined.

The function

$$\overline{I}(a) = \inf_{\substack{R: \int z_1 dR = a \\ R \text{ ergodic}}} H(R),$$

defined for $a \neq 0$, extends as a convex function to all of $R^d$ and, with this $\overline{I}$ as rate function, a large deviation result holds.

We now turn to the quenched case. Although a proof using the subadditive ergodic theorem exists, we will provide an alternate approach that is more appealing. We will illustrate this in the context of Brownian motion with a random drift. Instead of the action of $\mathbf{Z}^d$, we can have $\mathbf{R}^d$ acting on $(\Omega, \Sigma, P)$ ergodically and consider a diffusion on $\mathbf{R}^d$ with a random infinitesimal generator

$$(\mathcal{L}^\omega u)(x) = \tfrac{1}{2}(\Delta u)(x) + \langle b(\omega, x), (\nabla u)(x)\rangle$$

acting on smooth functions on $\mathbf{R}^d$. Here, $b(\omega, x)$ is generated from a map $b(\omega): \Omega \to \mathbf{R}^d$ by the action of $\{\tau_x : x \in \mathbf{R}^d\}$ on $\omega$:

$$b(\omega, x) = b(\tau_x \omega).$$

Again, there is the quenched measure $Q^\omega$ that corresponds to the diffusion with generator $\mathcal{L}^\omega$ that starts from 0 at time 0, and the averaged measure that is given by a similar formula. This model is referred to as *diffusion with a random drift*. Exactly the same questions can be asked in this context. We can define a diffusion on $\Omega$ with generator

$$\mathcal{L} = \tfrac{1}{2}\Delta + \langle b(\omega), \nabla\rangle,$$

where $\nabla = \{D_i\}$ are the generators of the translation group $\{\tau_x : x \in \mathbf{R}^d\}$. This is essentially the image of lifting the paths $x(t)$ of the diffusion on $\mathbf{R}^d$ corresponding to $\mathcal{L}^\omega$ to $\Omega$ by

$$\omega(t) = \tau_{x(t)}\omega.$$

While there is no possibility of having an invariant probability measure on $\mathbf{R}^d$, on $\Omega$, one can hope to find an invariant probability density $\phi(\omega)$, that is, to find $\phi(\omega) \geq 0$ in $L_1(P)$ with $\int \phi \, dP = 1$ which solves

$$\tfrac{1}{2}\Delta\phi = \nabla \cdot (b\phi).$$

If such a $\phi$ exists, then we have an ergodic theorem for the diffusion process $Q^\omega$ corresponding to $\mathcal{L}$ on $\Omega$,

$$\lim_{t\to\infty} \frac{1}{t} \int_0^t f(\omega(s)) \, ds = \int f(\omega)\phi(\omega) \, dP \qquad \text{a.e. } P.$$

This also translates to an ergodic theorem on $R^d$. If we define the stationary process $g$ by $g(\omega, x) = f(\tau_x \omega)$, then

$$\lim_{t\to\infty} \frac{1}{t} \int_0^t f(\omega, x(s)) \, ds = \int f(\omega)\phi(\omega) \, dP \qquad \text{a.e. } Q^\omega, \text{ a.e. } P,$$

where $Q^\omega$ is now the quenched process in the random environment. Since

$$x(t) = \int_0^t b(\omega, x(s)) \, ds + \beta(t),$$

where $\beta(\cdot)$ is the Brownian motion, it is clear that

$$\lim_{t\to\infty} \frac{x(t)}{t} = \int b(\omega)\phi(\omega) \, dP \qquad \text{a.e. } Q^\omega, \text{ a.e. } P,$$

providing a law of large numbers for $x(t)$. While we cannot be sure of finding $\phi$ for a given $b$, it is easy to find a $b$ for a given $\phi$. For instance, if $\phi > 0$, we could take $b = \frac{\nabla\phi}{2\phi}$. Or, more generally, $b = \frac{\nabla\phi}{2\phi} + \frac{c}{\phi}$ with $\nabla \cdot c = 0$. If we change $b$ to $b'$ which satisfies $\frac{1}{2}\Delta\phi = \nabla \cdot (b'\phi)$, the new process $Q^{b',\omega}$ with drift $b'$ will, in the time interval $[0, t]$, have relative entropy

$$E^{Q^{b',\omega}}\left[\tfrac{1}{2} \int_0^t \|b(\omega(s)) - b'(\omega(s))\|^2 \, ds\right]$$

and, by the ergodic theorem, one can see that, a.e. with respect to $P$,

$$\lim_{t\to\infty} \frac{1}{t} E^{Q^{b',\omega}}\left[\frac{1}{2} \int_0^t \|b(\omega(s)) - b'(\omega(s))\|^2 \, ds\right] = \frac{1}{2} \int \|b(\omega) - b'(\omega)\|^2 \phi(\omega) \, dP.$$

Moreover, for almost all $\omega$ with respect to $P$, almost surely with respect to $Q^{b',\omega}$,

$$\lim_{t\to\infty} \frac{x(t)}{t} = \int b'(\omega)\phi(\omega) \, dP.$$

If we fix $\int b'(\omega)\phi(\omega) \, dP = a$, the bound

$$\liminf_{t\to\infty} \frac{1}{t} \log Q^\omega\left[\frac{x(t)}{t} \simeq a\right] \geq -\frac{1}{2} \int \|b - b'\|^2 \phi \, dP$$

is easily obtained. If we define

$$I(a) = \inf_{\substack{b',\phi \\ (1/2)\Delta\phi = \nabla \cdot (b'\phi) \\ \int b'\phi\, dP = a}} \tfrac{1}{2} \int \|b - b'\|^2 \phi\, dP,$$

then

$$\liminf_{t \to \infty} \frac{1}{t} \log Q^\omega \left[ \frac{x(t)}{t} \simeq a \right] \geq -I(a).$$

Of course, these statements are valid a.e. $P$. One can check that $I$ is convex and that the upper bound amounts to proving the dual estimate

$$\lim_{t \to \infty} \frac{1}{t} \log E^{Q^\omega}[e^{\langle \theta, x(t) \rangle}] \leq \Psi(\theta),$$

where

$$\Psi(\theta) = \sup_a [\langle a, \theta \rangle - I(a)].$$

We need a bound on the solution of

$$u_t = \tfrac{1}{2}\Delta u + \langle b, \nabla u \rangle$$

with $u(0) = \exp[\langle \theta, x \rangle]$. By the Hopf–Cole transformation $v = \log u$, this reduces to estimating

$$v_t = \tfrac{1}{2}\Delta v + \tfrac{1}{2}\|\nabla v\|^2 + \langle b, \nabla v \rangle$$

with $v(0) = \langle \theta, x \rangle$. This can be done if we can construct a subsolution

$$\tfrac{1}{2}\nabla \cdot w + \tfrac{1}{2}\|\nabla w\|^2 + \langle b, w \rangle \leq \psi(\theta)$$

on $\Omega$, where $w : \Omega \to \mathbf{R}^d$ satisfies $\int w\, dP = \theta$ and $\nabla \times w = 0$ in the sense that $D_i w_j = D_j w_i$. The existence of the subsolution comes from convex analysis.

$$\psi(\theta) = \sup_{\substack{b',\phi \\ (1/2)\Delta\phi = \nabla \cdot (b'\phi)}} \left[ \int \langle b', \theta \rangle \phi\, dP - \tfrac{1}{2} \int \|b - b'\|^2 \phi\, dP \right]$$

$$= \sup_\phi \sup_{b'} \inf_u \left[ \int \langle b', \theta \rangle \phi\, dP - \tfrac{1}{2} \int \|b - b'\|^2 \phi\, dP \right.$$

$$\left. + \tfrac{1}{2}[\Delta u + \langle b', \nabla u \rangle]\phi\, dP \right]$$

$$= \sup_\phi \inf_u \sup_{b'} \left[ \int \langle b', \theta \rangle \phi\, dP - \tfrac{1}{2} \int \|b - b'\|^2 \phi\, dP \right.$$

$$\left. + \int \tfrac{1}{2}[\Delta u + \langle b', \nabla u \rangle]\phi\, dP \right]$$

$$= \sup_\phi \inf_u \int \sup_{b'} \left[ \langle b', \theta + \nabla \cdot u \rangle \phi \, dP - \tfrac{1}{2} \int \|b - b'\|^2 \phi \, dP + \int \tfrac{1}{2} \Delta u \phi \, dP \right]$$

$$= \sup_\phi \inf_u \int [\tfrac{1}{2}\Delta u + \langle b, \theta + \nabla u\rangle + \tfrac{1}{2}\|\theta + \nabla u\|^2] \phi \, dP$$

$$= \sup_\phi \inf_{\substack{\nabla \times w = 0 \\ \int w \, dP = \theta}} \int [\tfrac{1}{2} \nabla \cdot w + \langle b, w \rangle + \tfrac{1}{2}\|w\|^2] \phi \, dP$$

$$= \inf_{\substack{\nabla \times w = 0 \\ \int w \, dP = \theta}} \sup_\phi \int [\tfrac{1}{2} \nabla \cdot w + \langle b, w \rangle + \tfrac{1}{2}\|w\|^2] \phi \, dP$$

$$= \inf_{\substack{\nabla \times w = 0 \\ \int w \, dP = \theta}} \sup_\omega [\tfrac{1}{2} \nabla \cdot w + \langle b, w \rangle + \tfrac{1}{2}\|w\|^2],$$

which proves the existence of a subsolution. One needs to justify the free interchange of sup and inf. In passing from the first line to the second, the restriction on $b'$ and $\phi$ is replaced by the Lagrange multiplier $u$.

This can be viewed as showing the existence of a limit as $\epsilon \to 0$ (homogenization) of the solution of

$$u_t^\epsilon = \frac{\epsilon}{2} \Delta u^\epsilon + \frac{1}{2} \|\nabla u^\epsilon\|^2 + \left\langle b\left(\frac{x}{\epsilon}, \omega\right), \nabla u^\epsilon \right\rangle$$

with $u^\epsilon(0, x) = f(x)$. The limit satisfies

$$u_t = \Psi(\nabla u)$$

with $u(0, x) = f(x)$.

**11. Homogenization of Hamilton–Jacobi–Bellman equations.** This can be generalized to equations of the form

$$u_t^\epsilon = \frac{\epsilon}{2} \Delta u^\epsilon + H\left(\frac{x}{\epsilon}, \nabla u^\epsilon, \omega\right), \qquad u(0,x) = f(x),$$

where $H(x, p, \omega) = H(\tau_x \omega, p)$ is a stationary process of convex functions in $p$. By the changes of variables $x = \epsilon y$, $u = \epsilon v$, $t = \epsilon \tau$, this reduces to the behavior of $\epsilon v^\epsilon(\frac{T}{\epsilon}, \frac{x}{\epsilon})$, where $v^\epsilon$ solves

$$v_\tau^\epsilon = \tfrac{1}{2} \Delta v^\epsilon + H(y, \nabla v^\epsilon, \omega), \qquad v^\epsilon(0, y) = \epsilon^{-1} f(\epsilon y).$$

From the principle of dynamic programming, we can represent the solution as the supremum of a family of solutions of linear equations,

$$v^\epsilon(\tau, y) = \sup_{b(\cdot, \cdot) \in \mathcal{B}} w^\epsilon(b(\cdot, \cdot), T - \tau, y),$$

where $w^\epsilon$ solves

$$w^\epsilon_\tau + \tfrac{1}{2}\Delta w^\epsilon + \langle b(\tau,y), \nabla w^\epsilon \rangle - L(\tau_y\omega, b(\tau,y)) = 0,$$
$$w^\epsilon(\epsilon^{-1}T, y) = \epsilon^{-1} f(\epsilon y),$$

$L$ being the convex dual of $H$ with respect to the variable $p$. If $b$ is chosen as $b(\tau_y\omega)$, for suitable choices of $b(\omega) \in \mathcal{C}$ that admit positive integrable solutions $\phi$ to

$$\tfrac{1}{2}\Delta\phi(\tau_x\omega) = \nabla \cdot b(\tau_x\omega)\phi(\tau_x,\omega),$$

then it is not hard to see that

$$\epsilon w^\epsilon(0,0) \to f\left(T \int b(\omega)\phi(\omega)\,dP\right) - T \int L(b(\omega),\omega)\phi(\omega)\,dP.$$

This provides a lower bound

$$\liminf_{\epsilon \to 0} u^\epsilon(T,0) \geq \sup_{b \in \mathcal{C}} \left[ f\left(T \int b(\omega)\phi(\omega)\,dP\right) - T \int L(b(\omega),\omega)\phi(\omega)\,dP \right]$$

which can be shown to be an upper bound as well.

**12. History and references.** The origin of large deviation theory goes back to Scandinavian actuaries [10] who were interested in the analysis of risk in the insurance industry. For sums of independent random variables, the general large deviations result was established by Cramér in [1]. The result for empirical distributions of independent identically distributed random variables is due to Sanov [18]. The generalization to Markov chains and processes can be found in several papers of Donsker and Varadhan [3, 4, 5, 6] and Gärtner [13]. The results concerning small random perturbations of deterministic systems goes back to the work of Varadhan [20], as well as Vencel and Freidlin [12]. Several monographs have appeared on the subject. Lecture notes by Varadhan [21], texts by Deuschel and Stroock [7], Dembo and Zeitouni [2], Schwartz and Weiss [19], Ellis [9], Dupuis and Ellis [8] and, most recently, by Feng and Kurtz [11]. They cover a wide spectrum of topics in large deviation theory. For large deviations in the context of hydrodynamic scaling, there is the text by Kipnis and Landim [14], as well as an exposition by Varadhan [23]. As for large deviations for random walks in a random environment, see [24], as well as references in Zeitouni's article [25]. For a general survey on large deviations and entropy, see [22]. The results on homogenization of random Hamilton–Jacobi–Bellman equations and its application to large deviations has appeared in [15] and [17]. Undoubtedly, there are many more references. A recent Google search on "Large Deviations" produced 3.4 million hits.

COURANT INSTITUTE OF MATHEMATICAL SCIENCES
NEW YORK UNIVERSITY
NEW YORK, NEW YORK 10012
USA
E-MAIL: varadhan@cims.nyu.edu